\documentclass[10pt]{article}

\usepackage{amsmath,amsbsy,amssymb,latexsym,amsthm,dsfont}
\usepackage{pstricks}
\usepackage{paralist}
\usepackage{graphics}
\usepackage{epsfig}

\usepackage{cite}

\newtheorem{theorem}{Theorem}[section]

\theoremstyle{definition}
\newtheorem{definition}[theorem]{Definition}
\newtheorem{remark}{Remark}
\newtheorem{example}{Example}

\newenvironment{keywords}{\begin{center}
\begin{minipage}[c]{11.2cm}
{\bf Keywords:}} {\end{minipage}
\end{center}}

\newenvironment{msc}{\begin{center}
\begin{minipage}[c]{11.2cm}
{\bf Mathematics Subject Classification 2010:}}
{\end{minipage}
\end{center}}

\title{Necessary and sufficient conditions\\
for the fractional calculus of variations\\
with Caputo derivatives\footnote{Submitted 6/March/2010
to \emph{Communications in Nonlinear Science and Numerical Simulation}; 
revised 12/July/2010; accepted for publication 16/July/2010.}}

\author{Ricardo Almeida\\
\texttt{ricardo.almeida@ua.pt}
\and Delfim F. M. Torres\\
\texttt{delfim@ua.pt}}

\date{Department of Mathematics\\
University of Aveiro\\
3810-193 Aveiro, Portugal}

\begin{document}

\maketitle

\begin{abstract}
We prove optimality conditions for different variational functionals
containing left and right Caputo fractional derivatives.
A sufficient condition of minimization under
an appropriate convexity assumption is given. An Euler--Lagrange equation
for functionals where the lower and upper bounds of the integral
are distinct of the bounds of the Caputo derivative is also proved.
Then, the fractional isoperimetric problem is formulated with
an integral constraint also containing Caputo derivatives.
Normal and abnormal extremals are considered.
\end{abstract}

\begin{msc}
49K05, 26A33.
\end{msc}



\section{Introduction}

Fractional Calculus is the branch of Mathematics that generalizes the derivative
and the integral of a function to a non-integer order. The subject is not recent
and it is as old as the calculus itself. In a letter dated 30th September 1695,
Leibniz proposed the following problem to L'Hopital: "Can the meaning of derivatives
with integer order be generalized to derivatives with non-integer orders?" Since then,
several mathematicians studied this question, among them Liouville, Riemann, Weyl and Letnikov.
There are many fields of applications where we can use the fractional calculus, like viscoelasticity,
electrochemistry, diffusion processes, control theory, heat conduction, electricity, mechanics,
chaos and fractals (see some references at the end, \textrm{e.g.},
\cite{Blutzer,Jonsson,Lorenzo,Mainardi,Miller,Podlubny,Sabatier,samko}).
To solve fractional differential equations, there exist several methods:
Laplace and Fourier transforms, truncated Taylor series, numerical methods, etc.
(see \cite{AGRA4} and references therein).
Recently, a lot of attention has been put on the fractional calculus of variations
(see, \textrm{e.g.}, \cite{AGRA1,AGRA2,AGRA3,Almeida2,Atanackovic,Dreisigmeyer1,%
Dreisigmeyer2,Frederico:Torres07,Klimek1,Klimek2,Muslih1,Muslih2,Riewe1,Riewe2}.
We also mention \cite{Almeida1}, were necessary and sufficient conditions of optimality
for functionals containing fractional integrals and fractional derivatives are presented.
For results on fractional optimal control see \cite{AGRA0,Frederico:Torres08}.
In the present paper we work with the Caputo fractional derivative. For problems of calculus of variations
with boundary conditions, Caputo's derivative seems to be more natural,
since for a given function $y$ to have continuous Riemann--Liouville fractional derivative
on a closed interval $[a,b]$, the function must satisfy the conditions $y(a)=y(b)=0$ \cite{Almeida2}.
We also mention that, if $y(a)=0$, then the left Riemann--Liouville derivative
is equal to the left Caputo derivative.

The paper is organized as follows. In Section~\ref{sec2} we present the necessary definitions
and some necessary facts about fractional calculus.
Section~\ref{sec3} is dedicated to our original results. We study fractional
Euler--Lagrange equations and the fractional isoperimetric problem within Caputo's
fractional derivative context for different kinds of functionals. We also give
sufficient conditions of optimality for fractional variational problems.


\section{Preliminaries}
\label{sec2}

\subsection{Review on fractional calculus}

There are several definitions of fractional derivatives and fractional integrals,
like Riemann--Liouville, Caputo, Riesz, Riesz–Caputo, Weyl, Grunwald--Letnikov, Hadamard, Chen, etc.
We will present the definitions of the first two of them. Except otherwise stated,
proofs of results may be found in \cite{kilbas}.

Let $f:[a,b]\rightarrow\mathbb{R}$ be a function, $\alpha$ a positive real number,
$n$ the integer satisfying $n-1 \leq \alpha<n$, and $\Gamma$ the Euler gamma function. Then,
\begin{enumerate}
\item{the left and right Riemann--Liouville fractional integrals of order $\alpha$ are defined by
$${_aI_x^\alpha}f(x)=\frac{1}{\Gamma(\alpha)}\int_a^x (x-t)^{\alpha-1}f(t)dt,$$
and
$${_xI_b^\alpha}f(x)=\frac{1}{\Gamma(\alpha)}\int_x^b(t-x)^{\alpha-1} f(t)dt,$$
respectively;}
\item{the left and right Riemann--Liouville fractional derivatives of order $\alpha$ are defined by
$${_aD_x^\alpha}f(x)=\frac{d^n}{dx^n}  {_aI_x^{n-\alpha}}f(x) =\frac{1}{\Gamma(n-\alpha)}\frac{d^n}{dx^n}\int_a^x
(x-t)^{n-\alpha-1}f(t)dt,$$
and
$${_xD_b^\alpha}f(x)=(-1)^n\frac{d^n}{dx^n}  {_xI_b^{n-\alpha}}f(x)=\frac{(-1)^n}{\Gamma(n-\alpha)}\frac{d^n}{dx^n}\int_x^b
(t-x)^{n-\alpha-1} f(t)dt,$$
respectively;}
\item{the left and right Caputo fractional derivatives of order $\alpha$ are defined by
$${_a^CD_x^\alpha}f(x)= {_aI_x^{n-\alpha}}\frac{d^n}{dx^n} f(x)=\frac{1}{\Gamma(n-\alpha)}\int_a^x (x-t)^{n-\alpha-1}f^{(n)}(t)dt,$$
and
$${_x^CD_b^\alpha}f(x)=(-1)^n{_xI_b^{n-\alpha}}\frac{d^n}{dx^n}  f(x)=\frac{1}{\Gamma(n-\alpha)}\int_x^b(-1)^n
(t-x)^{n-\alpha-1} f^{(n)}(t)dt,$$
respectively.}
\end{enumerate}

There exists a relation between the  Riemann--Liouville and the  Caputo fractional derivatives:

$${_a^CD_x^\alpha}f(x)={_aD_x^\alpha}f(x)-\sum_{k=0}^{n-1}\frac{f^{(k)}(a)}{\Gamma(k-\alpha+1)}(x-a)^{k-\alpha}$$
and
$${_x^CD_b^\alpha}f(x)={_xD_b^\alpha}f(x)-\sum_{k=0}^{n-1}\frac{f^{(k)}(b)}{\Gamma(k-\alpha+1)}(b-x)^{k-\alpha}.$$
Therefore,
$$\mbox{if } f(a)=f'(a)=\ldots=f^{(n-1)}(a)=0, \mbox{ then } {_a^CD_x^\alpha}f(x)={_aD_x^\alpha}f(x)$$
and
$$\mbox{if } f(b)=f'(b)=\ldots=f^{(n-1)}(b)=0, \mbox{ then } {_x^CD_b^\alpha}f(x)={_xD_b^\alpha}f(x).$$
These fractional operators are linear, \textrm{i.e.},
$$\mathcal{P} (\mu f(x)+\nu g(x))=\mu\, \mathcal{P}f(x)+\nu \,\mathcal{P}g(x),$$
where $\mathcal{P}$ is ${_aD_x^\alpha}, {_xD_b^\alpha}, {_a^CD_x^\alpha}, {_x^CD_b^\alpha},{_aI_x^\alpha}$
or ${_xI_b^\alpha}$, and $\mu$ and $\nu$ are real numbers.

If $f \in C^n[a,b]$, then the left and right Caputo derivatives are continuous on $[a,b]$.
The main advantage of Caputo's approach is that the initial conditions for fractional
differential equations with Caputo derivatives take on the same form as for integer-order differential equations.

Some properties valid for integer differentiation and integer integration remain valid for fractional
differentiation and fractional integration; namely the Caputo fractional derivatives and the
Riemann--Liouville fractional integrals are inverse operations:
\begin{enumerate}
\item{If $f \in L_\infty(a,b)$ or $f \in C[a,b]$, and if $\alpha>0$, then
$${_a^CD_x^\alpha}{_aI_x^\alpha}f(x)=f(x) \quad \mbox{and} \quad {_x^CD_b^\alpha}{_xI_b^\alpha}f(x)=f(x).$$}
\item{If $f \in C^n[a,b]$ and if $\alpha>0$, then
$${_aI_x^\alpha}{_a^CD_x^\alpha}f(x)=f(x)-\sum_{k=0}^{n-1}\frac{f^{(k)}(a)}{k!}(x-a)^{k}$$
and
$${_xI_b^\alpha}{_x^CD_b^\alpha}f(x)=f(x)-\sum_{k=0}^{n-1}\frac{(-1)^kf^{(k)}(b)}{k!}(b-x)^{k}.$$}
\end{enumerate}

We also need for our purposes integration by parts formulas.
For $\alpha>0$, we have (\textrm{cf.} \cite{AGRA3})
$$
\begin{array}{ll}
\displaystyle\int_{a}^{b}g(x)\cdot {_a^C D_x^\alpha}f(x)dx&=\displaystyle\int_a^b f(x)\cdot {_x D_b^\alpha} g(x)dx \\
&+ \displaystyle\sum_{j=0}^{n-1}\left[{_xD_b^{\alpha+j-n}}g(x) \cdot {_xD_b^{n-1-j}} f(x)\right]_a^b
\end{array}
$$
and
$$
\begin{array}{ll}\displaystyle\int_{a}^{b}g(x)\cdot {_x^C D_b^\alpha}f(x)dx&=\displaystyle\int_a^b f(x)\cdot {_a D_x^\alpha} g(x)dx\\
&+ \displaystyle\sum_{j=0}^{n-1} \left[(-1)^{n+j}{_aD_x^{\alpha+j-n}}g(x) \cdot {_aD_x^{n-1-j}} f(x)\right]_a^b,
\end{array}
$$
where ${_aD_x^{k}}g(x)={_aI_x^{-k}}g(x)$ and  ${_xD_b^{k}}g(x)={_xI_b^{-k}}g(x)$ if $k<0$.

Therefore, if $0<\alpha<1$, we obtain
\begin{equation}\label{IP3}\int_{a}^{b}g(x)\cdot {_a^C D_x^\alpha}f(x)dx
=\int_a^b f(x)\cdot {_x D_b^\alpha} g(x)dx+\left[{_xI_b^{1-\alpha}}g(x) \cdot f(x)\right]_a^b\end{equation}
and
\begin{equation}
\label{IP4}
\int_{a}^{b}g(x)\cdot {_x^C D_b^\alpha}f(x)dx
=\int_a^b f(x)\cdot {_a D_x^\alpha} g(x)dx-\left[{_aI_x^{1-\alpha}}g(x) \cdot f(x)\right]_a^b.
\end{equation}
Moreover, if $f$ is a function such that $f(a)=f(b)=0$, we have simpler formulas:
\begin{equation}
\label{IP1}
\int_{a}^{b}  g(x)\cdot{_a^C D_x^\alpha}f(x)dx=\int_a^b f(x)\cdot {_x D_b^\alpha} g(x)dx
\end{equation}
and
\begin{equation}
\label{IP2}
\int_{a}^{b}  g(x)\cdot{_x^C D_b^\alpha}f(x)dx=\int_a^b f(x)\cdot {_a D_x^\alpha} g(x)dx.
\end{equation}

\begin{remark}
As $\alpha$ goes to $1$, expressions (\ref{IP3}) and (\ref{IP4}) reduce to the classical integration by parts formulas:
$$\int_{a}^{b}g(x)\cdot f'(x)dx=-\int_a^b f(x)\cdot g'(x)dx+\left[g(x) \cdot f(x)\right]_a^b$$
and
$$\int_{a}^{b}g(x)\cdot (-f'(x))dx=\int_a^b f(x)\cdot g'(x)dx-\left[g(x) \cdot f(x)\right]_a^b,$$
respectively.
\end{remark}

\begin{remark} Observe that the left member of equations (\ref{IP1}) and (\ref{IP2}) contains
a Caputo fractional derivative, while the other one contains a Riemann--Liouville
fractional derivative. However, since there exists a relation between the two derivatives,
we could present this formula with only one fractional derivative,
although in this case the resulting equation will contain some extra terms.
\end{remark}


\subsection{Fractional Euler--Lagrange equations}

From now on we fix $0<\alpha,\beta<1$. Also, to simplify, we denote
$$[y](x):=(x,y(x),\, {_a^C D_x^\alpha} y(x),\,{_x^C D_b^\beta} y(x)).$$
Let
$$
{_a^{\alpha}E_b^{\beta}}=\left\{ y : [a,b]\to\mathbb R \, | \, {_a^C D_x^\alpha}y
\mbox{ and } {_x^C D_b^\beta}y  \mbox{ exist and are continuous on} \, [a,b] \right\}.
$$

\begin{definition}
The space of variations $^CVar(a,b)$ for the Caputo derivatives is defined by
$$^CVar(a,b)= \left\{ h \in {_a^{\alpha}E_b^{\beta}} \, | \,  h(a)=h(b)=0\right\}.$$
\end{definition}

We now present first order necessary conditions of optimality for functionals,
defined on ${_a^{\alpha}E_b^{\beta}}$, of the type
\begin{equation}
\label{P0}
J(y)=\int_a^b L[y](x)dx.
\end{equation}
We assume that the map $(x,y,u,v)\to L(x,y,u,v)$ is a function of class $C^1$.
Denoting by $\partial_i L$ the partial derivative of $L$ with respect to the ith
variable, $i=1,2,3,4$, we also assume that $\partial_3 L$ has continuous right
Riemann--Liouville fractional derivative of order $\alpha$ and $\partial_4 L$
has continuous left Riemann--Liouville fractional derivative of order $\beta$.

\begin{definition}
We say that $y$ is a local minimizer (respectively local maximizer)
of $J$ if there exists a $\delta>0$ such that $J(y)\leq J(y_1)$
(respectively $J(y)\geq J(y_1)$) for all $y_1$ such that $\|y-y_1\|<\delta$.
\end{definition}

In \cite{AGRA1} Agrawal considers the problem of finding extremals for functionals
containing left and right Riemann--Liouville fractional derivatives of the form
$$J(y)=\int_a^b L(x,y(x),\,{_aD_x^\alpha}y(x),\, {_xD_b^\beta}y(x))dx$$
and he derived an Euler--Lagrange equation for an extremum $y$ of $J$,
subject to the boundary conditions $y(a)=y_a$, $y(b)=y_b$:
$$
\partial_2 L+{_xD_b^\alpha} \partial_3 L+{_aD_x^\beta} \partial_4 L=0,
\quad \mbox{for all } x \in [a,b].
$$

In \cite{Atanackovic} a new type of functional is studied, in case where the lower
bound of the integral do not coincide with the lower bound of the fractional derivative:
$${J^*}(y)=\int_A^B L(x,y(x),\,{_a D_x^\alpha} y(x))dx,$$
where $[A,B]\subset[a,b]$. We also mention \cite{AGRA2},
where an Euler--Lagrange equation and a transversality condition are given,
for functionals with left Caputo derivatives and a boundary condition on the initial point $x=a$.

\begin{theorem}[\cite{AGRA2}]
Let $J$ be the functional defined by
$$
J(y) =\int_a^b L(x,y(x),\,{_a^CD_x^\alpha} y(x))dx .
$$
Let $y$ be a local minimizer of $J$ satisfying the boundary condition $y(a)=y_0$.
Then, $y$ satisfies the following conditions:
\begin{equation}
\label{Agra1}
\partial_2 L+{_xD_b^\alpha} \partial_3 L=0
\end{equation}
and
\begin{equation}
\label{Agra2}
\left.{_x I _b^{1-\alpha}} \partial_3 L\right| _{x=b}=0.
\end{equation}
\end{theorem}

Note that the Euler--Lagrange equation (\ref{Agra1}) contains
a Riemann--Liouville fractional derivative, although $J$ has only Caputo's derivative.
Also, the transversality condition (\ref{Agra2}), in general, contains fractional derivative terms.
Thus, in order to solve a fractional variational problem, it may be required fractional boundary conditions.
This result is then proven for functionals with higher order fractional derivatives.
In \cite{AGRA6} functionals with the left and right Caputo fractional derivatives are considered.
We include here a short proof. For more on the subject we refer the reader to \cite{Cresson,Inizan,AM:frac1}.

\begin{theorem}[\cite{AGRA6}]
\label{ELtheorem}
Let $J$ be the functional as in (\ref{P0}) and $y$ a local minimizer of
$J$ satisfying the boundary conditions $y(a)=y_a$ and $y(b)=y_b$.
Then, $y$ satisfies the Euler--Lagrange equation
\begin{equation}
\label{ELequation}
\partial_2 L+{_xD_b^\alpha} \partial_3 L+{_aD_x^\beta} \partial_4 L=0.
\end{equation}
\end{theorem}

\begin{proof}
Let $\epsilon$ be a small real parameter and $\eta \in {^C Var(a,b)}$.
Consider a variation of $y$; say $y+\epsilon\eta$. Since the
Caputo derivative operators are linear, we have
$$
J(y+\epsilon\eta)=\int_a^b L(x,y+\epsilon\eta,\,{_a^CD_x^\alpha} y
+\epsilon {_a^CD_x^\alpha}\eta,\, {_x^CD_b^\beta}y+\epsilon {_x^CD_b^\beta}\eta)dx.
$$

We can regard $J$ as a function of one variable, $\hat J(\epsilon)=J(y+\epsilon\eta)$.
Since $y$ is the local minimizer, $\hat J$ attains an extremum at $\epsilon=0$.
Differentiating $\hat J(\epsilon)$ at zero, it follows that
$$\int_a^b \left[\partial_2 L \, \eta+ \partial_3 L \, {_a^CD_x^\alpha}\eta
+\partial_4 L \, {_x^CD_b^\beta}\eta \right]dx=0.$$
Integrating by parts, and since $\eta(a)=\eta(b)=0$, one finds that
$$
\int_a^b \left[ \partial_2 L + {_xD_b^\alpha} \partial_3 L
+ {_aD_x^\beta} \partial_4 L\right] \eta \, dx=0
$$
for all $\eta \in {^C Var(a,b)}$. By the arbitrariness of $\eta$ and by the
fundamental lemma of the calculus of variations
(see, \textrm{e.g.}, \cite[p.~32]{Brunt}), it follows that
$$\partial_2 L + {_xD_b^\alpha} \partial_3 L + {_aD_x^\beta} \partial_4 L=0.$$
\end{proof}

When the term ${_x^CD_b^\beta}y$ is not present in the function $L$,
then equation (\ref{ELequation}) reduces to a simpler one:
$$\partial_2 L+{_xD_b^\alpha} \partial_3 L=0.$$
Moreover, if we allow $\alpha=1$, and since in that case the right Riemann--Liouville
fractional derivative is equal to $-d/dx$, we obtain the classical Euler--Lagrange equation:
$$\partial_2 L-\frac{d}{dx} \partial_3 L=0.$$


\section{Main results}
\label{sec3}

Our first contribution is to generalize Theorem~\ref{ELtheorem}
(\textrm{cf.} Theorem~\ref{E-L3} below). We consider a new type
of functional, where the lower bound of the integral is greater
than the lower bound of the Caputo's derivative, and the upper bound
of the integral is less than the upper bound of the Caputo's derivative.
Because of this, we can not apply directly the integration by parts formula
and some technical auxiliary procedures are required. A similar problem
is addressed for the Riemann--Liouville fractional derivatives in \cite{Atanackovic}.


\subsection{The Euler--Lagrange equation}

We consider the functional
\begin{equation}
\label{Funct}
{J^*}(y) =\int_A^B L(x,y(x),\,{_a^CD_x^\alpha} y(x),\, {_x^CD_b^\beta}y(x))dx,
\end{equation}
where $(x,y,u,v)\to L(x,y,u,v)\in C^1$ and $[A,B]\subset[a,b]$.

\begin{theorem}
\label{E-L3} If $y$ is a local minimizer of ${J^*}$ given by (\ref{Funct}),
satisfying the boundary conditions $y(a)=y_a$ and $y(b)=y_b$, then $y$ satisfies the system
$$
\left\{\begin{array}{ll}
\displaystyle \partial_2 L + {_xD_B^\alpha} \partial_3 L
+{_AD_x^\beta} \partial_4 L=0 & \mbox{ for all } x \in [A,B]\vspace{0.1cm}\\
\displaystyle {_xD_B^\alpha}\partial_3 L -{_xD_A^\alpha}\partial_3 L=0 & \mbox{ for all } x \in [a,A]\vspace{0.2cm}\\
\displaystyle {_AD_x^\beta} \partial_4 L -{_BD_x^\beta}\partial_4 L=0 & \mbox{ for all } x \in [B,b]
\end{array}.\right.
$$
\end{theorem}

\begin{proof}
Let $y$ be a minimizer and let $\hat y=y+\epsilon\eta$ be a variation of $y$,
$\eta \in {^C Var(a,b)}$, such that $\eta(A)=\eta(B)=0$. Define the new function
$\hat {J^*} (\epsilon)={J^*}(y+\epsilon\eta)$. By hypothesis, $y$ is a local extremum
of ${J^*}$ and so $\hat {J^*}$ has a local extremum at $\epsilon=0$. Therefore, the following holds:
$$
\begin{array}{ll}0&=\displaystyle \int_A^B \left[\partial_2 L \, \eta
+\partial_3 L\, {_a^CD_x^\alpha} \eta +\partial_4 L\, {_x^CD_b^\beta}\eta\right]dx\\
&=\displaystyle  \int_A^B \partial_2 L \, \eta dx
+\left[\int_a^B \partial_3 L\, {_a^CD_x^\alpha} \eta dx
- \int_a^A \partial_3 L\, {_a^CD_x^\alpha} \eta dx \right] \\
&\quad +\displaystyle  \left[ \int_A^b \partial_4 L\, {_x^CD_b^\beta}\eta dx
- \int_B^b \partial_4 L\, {_x^CD_b^\beta}\eta dx  \right].
\end{array}
$$
Integrating by parts the four last terms gives:
$$\begin{array}{ll}0&= \displaystyle \int_A^B \partial_2 L\, \eta dx
+\left[\int_a^B \eta \, {_xD_B^\alpha}\partial_3 Ldx - \int_a^A \eta \, {_xD_A^\alpha}\partial_3 Ldx \right]\\
&\quad + \displaystyle \left[ \int_A^b \eta \, {_AD_x^\beta}\partial_4 Ldx
- \int_B^b\eta \, {_BD_x^\beta}\partial_4 Ldx  \right]\\
&= \displaystyle \int_A^B \partial_2 L\, \eta dx +\left[\int_a^A \eta \, {_xD_B^\alpha}\partial_3 Ldx
+ \int_A^B \eta \, {_xD_B^\alpha}\partial_3 Ldx - \int_a^A \eta \, {_xD_A^\alpha}\partial_3 Ldx \right] \\
&\quad + \displaystyle \left[ \int_A^B \eta \, {_AD_x^\beta}\partial_4 Ldx
+\int_B^b \eta \, {_AD_x^\beta}\partial_4 Ldx - \int_B^b\eta \, {_BD_x^\beta}\partial_4Ldx  \right]\\
&= \displaystyle \int_a^A\left[ {_xD_B^\alpha}\partial_3 L- {_xD_A^\alpha}\partial_3 L\right] \eta dx
+ \int_A^B \left[\partial_2 L+{_xD_B^\alpha}\partial_3 L+ {_AD_x^\beta}\partial_4 L\right]\eta dx\\
&\quad +\displaystyle \int_B^b \left[{_AD_x^\beta}\partial_4 L- {_BD_x^\beta}\partial_4 L\right]\eta dx.
\end{array}
$$
Since $\eta$ is an arbitrary function, we can assume that $\eta(x)=0$ for all
$x \in [A,b]$ and so by the fundamental lemma of the calculus of variations,
$$
{_xD_B^\alpha}\partial_3 L- {_xD_A^\alpha}\partial_3 L=0,\mbox{ for all } x\in[a,A].
$$
Similarly, one proves the other two conditions:
$$
\partial_2 L+{_xD_B^\alpha}\partial_3 L+ {_AD_x^\beta}\partial_4 L=0,\mbox{ for all } x\in[A,B],
$$
and
$$
{_AD_x^\beta}\partial_4 L- {_BD_x^\beta}\partial_4 L=0,\mbox{ for all } x\in[B,b].
$$
\end{proof}


\subsection{The fractional isoperimetric problem}

The \emph{isoperimetric problem} is one of the most ancient problems of the calculus of variations.
For example, given a positive real number $l$, what is the shape of the closed curve $C$
of length $l$ which defines the maximal area? The most essential contribution towards
its rigorous proof was given in 1841 and is due to Jacob Steiner (1796--1863).
We state the fractional isoperimetric problem as follows.

Given a functional $J$ as in (\ref{P0}), which functions $y$
minimize (or maximize) $J$, when subject to given boundary conditions
\begin{equation}
\label{bouncond} y(a)=y_a,\ y(b)=y_b,
\end{equation}
and an integral constraint
\begin{equation}
\label{isocons} I(y)=\int_a^b g[y]dx=l.
\end{equation}
Here, similarly as before, we consider a function $g$ of class $C^1$, such that $\partial_3 g$
has continuous right Riemann--Liouville fractional derivative of order $\alpha$ and $\partial_4 g$
has continuous left Riemann--Liouville fractional derivative of order $\beta$, and functions
$y \in {_a^{\alpha}E_b^{\beta}}$. A function $y \in {_a^{\alpha}E_b^{\beta}}$ that satisfies
(\ref{bouncond}) and (\ref{isocons}) is called admissible.

\begin{definition}
\label{extermal}
An admissible function $y$ is an \emph{extremal} for $I$ in (\ref{isocons})
if it satisfies the equation
$$
\partial_2 g [y](x) +{_xD_b^\alpha} \partial_3 g [y](x)
+{_aD_x^\beta} \partial_4 g [y](x)=0,\ \mbox{for all } x\in[a,b].
$$
\end{definition}

Observe that Definition~\ref{extermal} makes sense by Theorem~\ref{ELtheorem}.
To solve the isoperimetric problem, the idea is to consider a new extended function.
The exact formula for such extended function depends on $y$ being or not an extremal
for the integral functional $I(y)$ (\textrm{cf.} Theorems~\ref{IsoPro} and \ref{IsoPro2}).

\begin{theorem}
\label{IsoPro}
Let $y$ be a local minimum for $J$ given by (\ref{P0}), subject to the conditions
(\ref{bouncond}) and (\ref{isocons}). If $y$ is not an extremal for the functional $I$,
then there exists a constant $\lambda$ such that $y$ satisfies
\begin{equation}
\label{E-L}
\partial_2 F+{_xD_b^\alpha} \partial_3 F+{_aD_x^\beta}\partial_4 F=0
\end{equation}
for all $x \in [a,b]$, where $F=L+\lambda g$.
\end{theorem}

\begin{proof}
Let $\eta_1,\eta_2 \in {^CVar(a,b)}$ be two functions, $\epsilon_1$ and $\epsilon_2$
two reals, and consider the new function of two parameters
\begin{equation}
\label{variation}
\hat{y}=y+\epsilon_1\eta_1+\epsilon_2\eta_2.
\end{equation}
The reason why we consider two parameters is because we can choose one of them
as a function of the other in order to $\hat y$ satisfy the integral constraint. Let
$$\hat{I}(\epsilon_1,\epsilon_2)=\int_a^b g[y+\epsilon_1\eta_1+\epsilon_2\eta_2](x)dx-l.$$
It follows by integration by parts that
\begin{align*}
\left.\frac{\partial \hat I}{\partial \epsilon_2} \right|_{(0,0)}
&=\int_a^b\left( \partial_2 g \, \eta_2+\partial_3 g \, {_a^CD_x^\alpha} \eta_2
+\partial_4 g \, {_x^CD_b^\beta} \eta_2 \right)dx\nonumber\\
&=\int_a^b \left(\partial_2 g+{_xD_b^\alpha}\partial_3 g+{_aD_x^\beta}\partial_4 g \right)\eta_2 dx.
\end{align*}
We have assumed that $y$ is not an extremal for $I$, and therefore there exists
a function $\eta_2$ satisfying the condition
\begin{equation}
\label{implicit}
\left.\frac{\partial \hat I}{\partial \epsilon_2} \right|_{(0,0)}\neq 0.
\end{equation}
Using (\ref{implicit}) and the fact that $\hat I(0,0)=0$,
by the Implicit Function Theorem there exists a  $C^1$ function $\epsilon_2(\cdot)$,
defined in some neighborhood of zero, such that
$$\hat I(\epsilon_1,\epsilon_2(\epsilon_1))=0.$$
Therefore, there exists a family of variations
of type (\ref{variation}) which satisfy the integral constraint.

We will now prove condition (\ref{E-L}). Similarly as before,
we define a new function of two variables $\hat J(\epsilon_1,\epsilon_2)=J(\hat{y})$.
Since $(0,0)$ is a local minimum of $\hat J$, subject to the constraint
$\hat I(0,0)=0$, and $\nabla \hat I(0,0)\neq \textbf{0}$,
by the Lagrange Multiplier Rule (see, \textrm{e.g.}, \cite[p.~77]{Brunt}),
there exists a constant $\lambda$ for which the following holds:
$$\nabla(\hat J(0,0)+\lambda \hat I(0,0))=\textbf{0}.$$

Simple calculations show that
$$
\left.\frac{\partial \hat J}{\partial \epsilon_1} \right|_{(0,0)}
=\int_a^b \left[ \partial_2 L+{_xD_b^\alpha} \partial_3 L+ {_aD_x^\beta}\partial_4 L\right]\eta_1 dx
$$
and
$$
\left.\frac{\partial \hat I}{\partial \epsilon_1} \right|_{(0,0)}
=\int_a^b \left[\partial_2 g+{_xD_b^\alpha}\partial_3 g
+{_aD_x^\beta}\partial_4 g\right]\eta_1 dx.
$$
In conclusion, it follows that
$$
\int_a^b \left[\partial_2 L+{_xD_b^\alpha}\partial_3 L+{_aD_x^\beta}\partial_4 L
+\lambda \left(\partial_2 g+{_xD_b^\alpha}\partial_3 g+{_aD_x^\beta}\partial_4 g\right)\right]\eta_1 dx=0.
$$
By the arbitrariness of  $\eta_1$ and the fundamental lemma of calculus of variations, one must have
$$\partial_2 L+{_xD_b^\alpha}\partial_3 L+{_aD_x^\beta}\partial_4 L
+\lambda \left(\partial_2 g+{_xD_b^\alpha}\partial_3 g+{_aD_x^\beta}\partial_4 g\right)=0.$$
This is equivalent to
$$\partial_2 F+{_xD_b^\alpha}\partial_3 F+{_aD_x^\beta}\partial_4 F=0.$$
\end{proof}

\begin{example}
\label{example}
Let $\overline y(x)={E_{\alpha}}(x^\alpha)$, $x \in [0,1]$, where ${E_\alpha}$ is the Mittag--Leffler function:
$${E_\alpha}(x)=\sum_{k=0}^\infty \frac{x^k}{\Gamma(\alpha k+1)}, \quad x \in \mathbb R, \alpha>0.$$
When $\alpha=1$, the  Mittag--Leffler function is the exponencial function, $E_1(x)=e^x$.

The left Caputo fractional derivative of $\overline y$ is $\overline y$ (\textrm{cf.} \cite[p.~98]{kilbas}),
$${^C_0 D_x^\alpha} \overline y=\overline y.$$

Consider the following fractional variational problem:
\begin{equation}
\label{example2}
\left\{\begin{array}{l}
\displaystyle J(y)=\int_0^1({^C_0 D_x^\alpha} y)^2 \, dx \quad \rightarrow \quad \mbox{extr},\\
\displaystyle I(y)=\int_0^1 \overline y \, {^C_0 D_x^\alpha} y \, dx = \int_0^1 (\overline y)^2 \, dx,\\
\displaystyle y(0)=1 \quad \mbox{and} \quad y(1)= {E_\alpha}(1).
\end{array}\right.
\end{equation}

The augmented function is
$$F(x,y,{^C_0D_x^\alpha} y, {^C_x D _1^\beta}y,\lambda)
=({^C_0D_x^\alpha} y)^2+\lambda\overline y \, {^C_0 D_x^\alpha} y$$
and the fractional Euler--Lagrange equation is
$$\partial_2 F+{_xD_1^\alpha}\partial_3F+{_xD_1^\beta}\partial_4F=0$$
\textrm{i.e.},
$${_xD_1^\alpha}(2\,{^C_0D_x^\alpha} y+\lambda \overline y)=0.$$
A solution of this problem is $\lambda=-2$ and $y=\overline y$.

Observe that, as $\alpha\to1$, the variational problem (\ref{example2}) becomes
$$\int_0^1 y'^2 \, dx \quad \rightarrow \quad \mbox{extr},$$
$$\int_0^1 \overline y \, y' \, dx =\frac12(e^2-1),$$
$$y(0)=1 \quad \mbox{and} \quad y(1)=e,$$
and the Euler--Lagrange equation is
\begin{equation}
\label{DE}
\partial_2 F -\frac{d}{dx}\partial_3F=0 \Leftrightarrow -\frac{d}{dx}(2 y'-2 \overline y)=0,
\end{equation}
where $F=y'^2-2\overline y y'$. Also, for $\alpha=1$, $\overline y(x)=e^x$,
which is obviously a solution of the differential equation
(\ref{DE}) (\textrm{cf.} Figure~\ref{fig:1}).
\begin{figure}
\begin{center}
\pspicture(7,7)(0.5,0.5)
\rput(3,3){\includegraphics[scale=0.35]{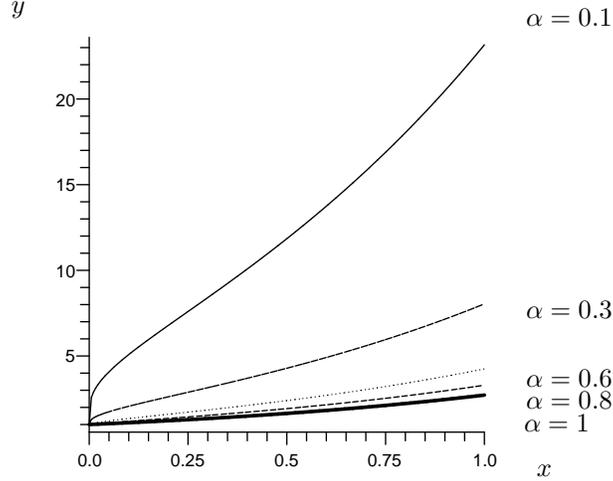}}
\uput{-1}[135](8.3,5.5){$\alpha=0.1$}
\uput{-1}[135](8.3,1.6){$\alpha=0.3$}
\uput{-1}[135](8.3,0.7){$\alpha=0.6$}
\uput{-1}[135](8.3,0.4){$\alpha=0.8$}
\uput{-1}[135](8,0.1){$\alpha=1$}
\uput{-1}[135](7.5,-0.5){$x$}
\uput{-1}[135](0.5,5.6){$y$}
\endpspicture
\end{center}
\caption{Solutions of problem (\ref{example2}).}
\label{fig:1}
\end{figure}
\end{example}

We now study the case when $y$ is an extremal of $I$
(the so called abnormal case).

\begin{theorem}
\label{IsoPro2}
Let $y$ be a local minimum of $J$ (\ref{P0}),
subject to the conditions (\ref{bouncond}) and (\ref{isocons}).
Then, there exist two constants $\lambda_0$ and $\lambda$, with $(\lambda_0,\lambda)\not=(0,0)$,
such that
$$\partial_2 K+{_xD_b^\alpha} \partial_3 K+{_aD_x^\beta}\partial_4 K=0$$
where $K=\lambda_0 L+\lambda g$.
\end{theorem}

\begin{proof}
Following the proof of Theorem \ref{IsoPro}, $(0,0)$ is an extremal of $\hat J$
subject to the constraint $\hat I=0$. Then, by the abnormal Lagrange multiplier rule
(see, \textrm{e.g.}, \cite[p.~82]{Brunt}), there exist two reals
$\lambda_0$ and $\lambda$, not both zero, such that
$$\nabla(\lambda_0 \hat J(0,0)+\lambda \hat I(0,0))=\textbf{0}.$$
Therefore,
$$\lambda_0 \displaystyle\left.\frac{\partial \hat J}{\partial \epsilon_1}\right|_{(0,0)}+
\lambda \displaystyle\left.\frac{\partial \hat I}{\partial \epsilon_1}\right|_{(0,0)}=0.$$
The rest of proof is similar to the one of Theorem \ref{IsoPro}.
\end{proof}


\subsection{An extension}

We now present a solution for the isoperimetric problem for functionals of type (\ref{Funct}).
Similarly, one has an integral constraint, but this time of the form
\begin{equation}
\label{Constraint}
{I^*}(y)=\int_A^B g(x,y(x),\, {_a^CD_x^\alpha} y(x),\,{_x^CD_b^\beta} y(x))dx=l.
\end{equation}

Again, we need the concept of extremal for a functional of type (\ref{Constraint}).

\begin{definition}
A function $y$ is called extremal for ${I^*}$ given by (\ref{Constraint}) if
$$\partial_2 g [y](x) +{_xD_B^\alpha} \partial_3 g[y](x)
+{_AD_x^\beta}\partial_4 g[y](x)=0,\ \mbox{for all } x\in[A,B].$$
\end{definition}

\begin{theorem}
\label{E-L5} If $y$ is a local minimum of ${J^*}$ given by (\ref{Funct}),
when restricted to the conditions  $y(a)=y_a$, $y(b)=y_b$ and
(\ref{Constraint}), and if $y$ is not an extremal for ${I^*}$,
then there exists a constant $\lambda$ such that
\begin{equation}
\label{IsoPro3}
\left\{\begin{array}{ll}
\displaystyle \partial_2 F+{_xD_B^\alpha} \partial_3 F+{_AD_x^\beta} \partial_4 F
=0 & \mbox{ for all } x \in [A,B]\vspace{0.1cm}\\
\displaystyle {_xD_B^\alpha} \partial_3 F -{_xD_A^\alpha} \partial_3 F=0 & \mbox{ for all } x \in [a,A]\vspace{0.2cm}\\
\displaystyle {_AD_x^\beta} \partial_4 F-{_BD_x^\beta} \partial_4 F=0 & \mbox{ for all } x \in [B,b]
\end{array}\right.
\end{equation}
where $F=L+\lambda g$.
\end{theorem}

\begin{remark}
In case $[A,B]=[a,b]$, Theorem~\ref{E-L5} is reduced to Theorem~\ref{IsoPro}.
\end{remark}

\begin{proof}
We consider a variation of form $\hat y =y+\epsilon_1 \eta_1+\epsilon_2\eta_2$, where
\begin{center}
$\eta_1,\eta_2 \in {^C Var(a,b)}$ and $\eta_1(A)=\eta_1(B)=\eta_2(A)=\eta_2(B)=0$.
\end{center}
Define $\hat {I^*}$ by the expression
$$\hat {I^*}(\epsilon_1,\epsilon_2)=\int_A^B g[\hat y]dx-l.$$

Then, $\hat {I^*}(0,0)=0$ and
$$
\begin{array}{ll}
\displaystyle \left.\frac{\partial \hat {I^*}}{\partial \epsilon_2} \right|_{(0,0)}
&=\displaystyle\int_A^B\left[\partial_2 g \, \eta_2+\partial_3 g \, {_a^CD_x^\alpha} \eta_2
+ \partial_4 g \, {_x^CD_b^\beta} \eta_2 \right]dx\\
&= \displaystyle \int_a^A\left[ {_xD_B^\alpha}\partial_3 g-{_xD_A^\alpha}\partial_3 g\right] \eta_2 dx
+ \int_A^B \left[\partial_2 g+{_xD_B^\alpha}\partial_3 g+ {_AD_x^\beta}\partial_4 g\right]\eta_2 dx\\
&\quad +\displaystyle \int_B^b \left[{_AD_x^\beta}\partial_4 g- {_BD_x^\beta}\partial_4 g\right]\eta_2 dx
\end{array}
$$
(the last expression follows by integration by parts and some technical calculations as presented
in the proof of Theorem~\ref{E-L3}). Let $\eta_2$ be a function such that
$$\left.\frac{\partial \hat {I^*}}{\partial \epsilon_2} \right|_{(0,0)}\neq 0$$
(its existence is guarantied since $y$ is not an extremal for ${I^*}$).
Therefore, we can consider a subset of the family of functions $\{ y+\epsilon_1 \eta_1
+ \epsilon_2\eta_2 \, | \, (\epsilon_1,\epsilon_2)\in \mathbb R^2 \}$
that is admissible for the isoperimetric problem. Let $\hat {J^*}(\epsilon_1,\epsilon_2)={J^*}(\hat{y})$.
Then, there exists a real $\lambda$ such that
\begin{equation}
\label{gradient}
\nabla(\hat {J^*}(0,0)
+\lambda \hat {I^*}(0,0))=\textbf{0}.
\end{equation}
Similarly, one has
$$
\begin{array}{ll}\displaystyle\left.\frac{\partial \hat {J^*}}{\partial \epsilon_1} \right|_{(0,0)}
&=\displaystyle \int_a^A\left[ {_xD_B^\alpha}\partial_3 L- {_xD_A^\alpha}\partial_3 L\right] \eta_1 dx
+ \int_A^B \left[\partial_2 L+{_xD_B^\alpha}\partial_3 L+{_AD_x^\beta}\partial_4 L\right]\eta_1 dx\\
 &\quad +\displaystyle\int_B^b \left[{_AD_x^\beta}\partial_4 L-{_BD_x^\beta}\partial_4 L\right]\eta_1 dx,
\end{array}
$$
and
$$
\begin{array}{ll}\displaystyle\left.\frac{\partial \hat {I^*}}{\partial \epsilon_1} \right|_{(0,0)}
&=\displaystyle \int_a^A\left[ {_xD_B^\alpha}\partial_3 g- {_xD_A^\alpha}\partial_3 g\right] \eta_1 dx
+ \int_A^B \left[\partial_2 g+{_xD_B^\alpha}\partial_3 g+ {_AD_x^\beta}\partial_4 g\right]\eta_1 dx\\
 &\quad +\displaystyle\int_B^b \left[{_AD_x^\beta}\partial_4 g-{_BD_x^\beta}\partial_4 g\right]\eta_1 dx.
\end{array}
$$
By (\ref{gradient}),
$$
\left.\frac{\partial \hat {J^*}}{\partial \epsilon_1} \right|_{(0,0)}
+\lambda \left.\frac{\partial \hat {I^*}}{\partial \epsilon_1} \right|_{(0,0)}=0
$$
and (\ref{IsoPro3}) follows from the arbitrariness of $\eta_1$.
\end{proof}

The following result generalizes Theorem \ref{IsoPro2} and is proved in a similar way.

\begin{theorem}
If $y$ is a local minimum of ${J^*}$ given by (\ref{Funct}), subject to the boundary conditions
(\ref{bouncond}) and the integral constraint (\ref{Constraint}),
then there exist two constants $\lambda_0$ and $\lambda$, not both zero, such that
$$
\left\{\begin{array}{ll}
\displaystyle \partial_2 K+{_xD_B^\alpha}\partial_3 K
+{_AD_x^\beta}\partial_4 K=0 & \mbox{ for all } x \in [A,B]\vspace{0.1cm}\\
\displaystyle{_xD_B^\alpha} \partial_3 K-{_xD_A^\alpha}\partial_3 K=0
& \mbox{ for all } x \in [a,A]\vspace{0.2cm}\\
\displaystyle{_AD_x^\beta} \partial_4 K-{_BD_x^\beta}\partial_4 K=0
& \mbox{ for all } x \in [B,b]\\
\end{array}\right.
$$
where $K=\lambda_0 L+\lambda g$.
\end{theorem}


\subsection{Sufficient conditions of optimality}

We are now interested in finding sufficient conditions for $J$ to attain local extremes.
Typically, some conditions of convexity over the Lagrangian are needed.

\begin{definition} We say that $f(\underline x,y,u,v)$
is convex in $S\subseteq\mathbb R^4$ if $\partial_2 f$, $\partial_3 f$ and $\partial_4 f$
exist and are continuous, and the condition
\begin{multline*}
f(x,y+y_1,u+u_1,v+v_1)-f(x,y,u,v)\\
\geq
\partial_2 f(x,y,u,v)y_1+\partial_3 f(x,y,u,v)u_1+\partial_4 f(x,y,u,v)v_1
\end{multline*}
holds for every $(x,y,u,v),(x,y+y_1,u+u_1,v+v_1)\in S$.
\end{definition}

\begin{theorem}
\label{SufficientCond}
Suppose that the function $L(\underline x,y,u,v)$ is convex in $[a,b]\times\mathbb R^3$.
Then each solution $y_0$ of the fractional Euler--Lagrange equation (\ref{ELequation}) minimizes
$$J(y)=\int_a^b L[y]dx,$$
when restricted to the boundary conditions $y(a)=y_0(a)$ and $y(b)=y_0(b)$.
\end{theorem}

\begin{proof}
Let $\eta\in{_a^{\alpha}E_b^{\beta}}$ be a function such that $\eta(a)=\eta(b)=0$.
Then, using integration by parts, it follows that
$$
\begin{array}{ll}
J(y_0+\eta)-J(y_0)&= \displaystyle\int_a^b(L[y_0+\eta]-L[y_0]) \, dx \\
& \geq \displaystyle \int_a^b \left( \partial_2 L[y_0]\eta
+ \partial_3 L[y_0] {_a^CD_x^\alpha}\eta +\partial_4 L[y_0]{_x^CD_b^\beta}\eta \right) \, dx\\
&=  \displaystyle\int_a^b \left[ \partial_2 L+{_xD_b^\alpha}\partial_3 L
+{_aD_x^\beta}\partial_4 L\right] [y_0] \cdot\eta\,dx=0
\end{array}
$$
since $y_0$ is a solution of (\ref{ELequation}).
We just proved that $J(y_0+\eta)\geq J(y_0)$.
\end{proof}

This procedure can be easily adapted for the isoperimetric problem.

\begin{theorem} Suppose that, for some constant $\lambda$, the functions
$L(\underline x,y,u,v)$ and $\lambda g (\underline x,y,u,v)$ are convex
in $[a,b]\times\mathbb R^3$. Let $F=L+\lambda g$.
Then each solution $y_0$ of the fractional equation
$$\partial_2 F+{_xD_b^\alpha}\partial_3 F+{_aD_x^\beta}\partial_4 F=0$$
minimizes
$$J(y)=\int_a^b L[y]dx,$$
under the constraints $y(a)=y_0(a)$ and $y(b)=y_0(b)$ and
$$I(y)=\int_a^b g[y]dx=l,\, l \in \mathbb R.$$
\end{theorem}

\begin{proof}
Let us prove that $y_0$ minimizes
$$\widetilde{F}(y)=\int_a^b(L[y]+ \lambda g[y] )\, dx.$$ First,
it is easy to prove that $L(\underline x,y,u,v)+\lambda g (\underline x,y,u,v)$ is convex.
Let $\eta \in {_a^{\alpha}E_b^{\beta}} $ be such that $\eta(a)=\eta(b)=0$.
Then, by Theorem~\ref{SufficientCond},
$\widetilde{F}(y_0+\eta)\geq \widetilde{F}(y_0)$. In other words,
if $y\in {_a^{\alpha}E_b^{\beta}}$ is any function
such that $y(a)=y_0(a)$ and $y(b)=y_0(b)$, then
$$
\int_a^b L[y]\, dx+\int_a^b \lambda g[y] \, dx
\geq \int_a^b L[y_0]\, dx+\int_a^b \lambda g[y_0] \, dx.
$$
If we restrict to the integral constraint, we obtain
$$\int_a^b L[y]\, dx+l \geq \int_a^b L[y_0]\, dx+l,$$
and so
$$\int_a^b L[y]\, dx\geq \int_a^b L[y_0]\, dx,$$
proving the desired result.
\end{proof}

\begin{example}
Recall Example~\ref{example}. Since $L(\underline x,y,u,v)=u^2$
and $\lambda g (\underline x,y,u,v)=-2\overline y u$ are both convex,
we conclude that $\overline y$ is actually a minimum
for the fractional variational problem (\ref{example2}).
\end{example}


\section*{Acknowledgements}

Work partially supported by the
{\it Centre for Research on Optimization and Control} (CEOC)
from the ``Funda\c{c}\~{a}o para a Ci\^{e}ncia e a Tecnologia'' (FCT),
cofinanced by the European Community Fund FEDER/POCI 2010.



\end{document}